\documentclass[smallextended]{svjour3}
\usepackage{amsmath,graphicx,indentfirst,enumitem,amssymb,
tikz,color,geometry,url}
\usepackage[subsection]{algorithm}
\usepackage{color}
\usepackage{amsfonts}
\usepackage{algorithm}
\usepackage{algpseudocode}
\usetikzlibrary{positioning,calc,arrows,automata}
\newcommand{\R}{\mathbb{R}}

\date{}
\begin{document}
\author{Silvia Noschese \and Lothar Reichel}
\institute{Silvia Noschese\at
Dipartimento di Matematica\\ 
SAPIENZA Universit\`a di Roma\\
P.le Aldo Moro 5, 00185 Roma, Italy\\
\email{noschese@mat.uniroma1.it}\\
Lothar Reichel\at
Department of Mathematical Sciences\\ 
Kent State University\\
Kent, OH 44242, USA\\
\email{reichel@math.kent.edu}
}
\title{Network analysis with the aid of the path length matrix}
\maketitle

\begin{abstract}
Let a network be represented by a simple graph $\mathcal{G}$ with $n$ vertices. A common 
approach to investigate properties of a network is to use the adjacency matrix 
$A=[a_{ij}]_{i,j=1}^n\in\R^{n\times n}$ associated with the graph $\mathcal{G}$, where 
$a_{ij}>0$ if there is an edge pointing from vertex $v_i$ to vertex $v_j$, and $a_{ij}=0$
otherwise. Both $A$ and its positive integer powers reveal important properties of the
graph. This paper proposes to study properties of a graph $\mathcal{G}$ by also using the 
path length matrix for the graph. The $(ij)^{th}$ entry of the path length matrix is 
the length of the shortest path from vertex $v_i$ to vertex $v_j$; if there is no path 
between these vertices, then the value of the entry is $\infty$. Powers of the path length 
matrix are formed by using min-plus matrix multiplication and are important for 
exhibiting properties of $\mathcal{G}$. We show how several known measures of 
communication such as closeness centrality, harmonic centrality, and eccentricity are 
related to the path length matrix, and we introduce new measures of communication, 
such as the harmonic $K$-centrality and global $K$-efficiency, where only (short)
paths made up of at most $K$ edges are taken into account. The sensitivity of the 
global $K$-efficiency to changes of the entries of the adjacency matrix also is 
considered.
\end{abstract}


\section{Introduction}\label{s1}
An important characteristic of a network is how well communication can flow in it, i.e.,
how easy or difficult it is to reach one part of the network from another part by 
following edges. How well information flows through the whole network can be measured by
the diameter of the graph that represents the network or by its global efficiency.
Both these measures can be determined with the aid of the path length matrix associated 
with the network. We will discuss these connections and introduce new measures of 
communication based on the path length matrix.

Let us introduce some notation and definitions that will be used throughout this paper. A 
network is represented by a graph. A \emph{weighted graph} 
$\mathcal{G}=\langle\mathcal{V},\mathcal{E},\mathcal{W}\rangle$ consists of a set of 
\emph{nodes} or \emph{vertices} $\mathcal{V}=\{v_1,v_2,\dots,v_n\}$, a set of \emph{edges} 
$\mathcal{E}=\{e_1,e_2,\dots,e_m\}$ that connect the vertices, and a set of weights 
$\mathcal{W}=\{a_{ij}\}_{i,j=1}^n$; the weights $a_{ij}$ are the entries of the 
adjacency matrix $A=[a_{ij}]_{i,j=1}^n$ associated with the graph $\mathcal{G}$; see 
below. An edge is said to be \emph{directed} if it starts at a vertex $v_i$ and ends
at a vertex $v_j$, and is denoted by $e(v_i\rightarrow v_j)$. An edge between the vertices 
$v_i$ and $v_j$ is said to be \emph{undirected} when the pair of vertices is unordered and
the weights $a_{ij}$ and $a_{ji}$ are positive and equal. An undirected egde between the 
vertices $v_i$ and $v_j$ is denoted by $e(v_i\leftrightarrow v_j)$. A graph with only 
undirected edges is said to be \emph{undirected}; otherwise the graph is \emph{directed}.
A \emph{simple graph} is a graph without multiple edges or self-loops. In particular, this
implies that the diagonal entries of the adjacency matrix for the graph vanish. This work 
considers simple graphs. 

The \emph{adjacency matrix} $A=[a_{ij}]_{i,j=1}^n\in\R^{n\times n}$ for a weighted graph 
$\mathcal{G}$ is determined by the weights $a_{ij}$ of the graph with $a_{ij}>0$ if there 
is an edge $e(v_i\rightarrow v_j)$ in $\mathcal{G}$. If there is no edge 
$e(v_i\rightarrow v_j)$ in $\mathcal{G}$, then $a_{ij}=0$. For an \emph{unweighted graph},
all positive entries $a_{ij}$ of $A$ equal one. A sequence of $k$ edges (not necessarily 
distinct) such that 
$\{e(v_1\rightarrow v_2),e(v_2\rightarrow v_3),\ldots,e(v_k\rightarrow v_{k+1})\}$ form a
\emph{walk}. If $v_{k+1}=v_1$, then the walk is said to be \emph{closed}. A sequence of 
distinct edges such that $\{e(v_1\rightarrow v_2),e(v_2\rightarrow v_3),\ldots,
e(v_k\rightarrow v_{k+1})\}$ form a \emph{path}. The \emph{length} of a path is given by 
the sum of all weights of the edges in the path.  (In the unweighted case, the sum of all 
weights of the edges in a path of length $k$ is $k$.) For further discussions on networks 
and graphs; see \cite{estrada2011structure,newman2010}.

To construct the path length matrix associated with the network, we will make use of 
{\it min-plus matrix multiplication}, i.e., matrix multiplication in the tropical algebra 
\cite{L}:
\begin{equation*}
C=A\star B: \qquad\qquad c_{ij}=\min_{h=1,2,\ldots,n} \{a_{ih}+b_{hj}\},\qquad 
1\leq i,j\leq n,
\end{equation*}
where $A,B,C\in\R^{n\times n}$.
We denote by $A^{1,\star}=[a^{(1,\star)}_{ij}]_{i,j=1}^n\in\R^{n\times n}$ the matrix 
obtained by setting to $\infty$ the vanishing off-diagonal entries of the adjacency 
matrix $A$ associated with the graph $\mathcal G$ under consideration. For $k>1$, the 
$k^{th}$ {\it min-plus power} of $A^{1,\star}$ is given by
\begin{equation*}
A^{k,\star}=[a^{(k,\star)}_{ij}]_{i,j=1}^n\in\R^{n\times n}: \qquad\qquad 
a^{(k,\star)}_{ij} =\min_{h=1,2,\ldots,n} \{a^{(k-1,\star)}_{ih}+a^{(1,\star)}_{hj}\}.
\end{equation*}
Notice that the matrix $A^{k,\star}$ gives vertex distances using paths of at most 
$k$ edges. In detail, the entry $a^{(k,\star)}_{ij}$, with $i\neq j$, represents the 
length of the \emph{shortest path} from $v_i$ to $v_j$ made up of at most $k$ edges. The 
diagonal entries of $A^{k,\star}$ are zero by definition. One has 
$a_{ij}^{(k,\star)}=\infty$ if every path from $v_i$ to $v_j$ is made up of more than $k$ 
edges, or if there is no path from $v_i$ to $v_j$.

The \emph{diameter} of a graph $\mathcal{G}$ is the maximal length $d_{\mathcal{G}}$ of 
the shortest path between any distinct vertices of the graph and provides a measure of how
easy it is for the vertices of the graph to communicate. One has
\begin{equation} \label{d_G}
d_{\mathcal{G}}=\max_{1\leq i,j\leq n} a_{ij}^{(n-1,\star)}.
\end{equation}
Indeed, the entry $a^{(n-1,\star)}_{ij}$ of the matrix 
$A^{n-1,\star}=[a^{(n-1,\star)}_{ij}]_{i,j=1}^n$ yields the length of the shortest path 
from $v_i$ to $v_j$. We will refer to $A^{n-1,\star}$ as the {\it path length matrix}.
Note that the triangle inequality 
holds for the entries of a path length matrix. Specifically, 
$$
a_{ij}^{(n-1,\star)}\leq
a_{ih}^{(n-1,\star)}+a_{hj}^{(n-1,\star)},\qquad 1\leq i,j\leq n.
$$

Consider an unweighted connected graph $\mathcal{G}$ with associated adjacency matrix 
$A\in\R^{n\times n}$. Recall that the diameter of $\mathcal{G}$ is the maximal number of 
edges in the shortest path between all pairs of distinct vertices of the graph. Given the 
vertices $v_i$ and $v_j$, there is an integer $\widehat{k}$, $1\leq\widehat{k}<n$, such that 
\[
a_{ij}^{(n-1,\star)}=\dots= a_{ij}^{(\widehat{k}+1,\star)}= a_{ij}^{(\widehat{k},\star)}=
\widehat{k}
\]
since the graph is connected, whereas for $1\leq h<\widehat{k}$, one has 
$a_{ij}^{(h,\star)}=\infty$. Thus, information provided by the path length matrix 
$A^{n-1,\star}$ includes information about all powers $A^{k,\star}$ for $1\leq k<n-1$.

Let $v_i$ and $v_j$ be distinct vertices in a weighted graph. Then there is an integer 
$\widehat{k}$, $1\leq\widehat{k}<n-1$, such that 
\[
a_{ij}^{(n-1,\star)}\leq \dots \leq a_{ij}^{(\widehat{k}+1,\star)}\leq 
a_{ij}^{(\widehat{k},\star)}<\infty
\]
and, for $1\leq h<\widehat{k}$, $a_{ij}^{(h,\star)}=\infty$. Thus, as in the unweighted 
case, information provided by the path length matrix $A^{n-1,\star}$ refines information 
given by the powers $A^{k,\star}$ for $1\leq k<n-1$. However, the information of the minimal 
number of steps required to reach vertex $v_j$ from vertex $v_i$ is lost.

As mentioned above, the path length matrix may be constructed by evaluating the min-plus 
powers of $A^{1,\star}$ $n-2$ times; here $A^{1,\star}$ is obtained from the adjacency 
matrix $A$ by setting all zero off-diagonal entries to $\infty$. The following MATLAB 
function, with the adjacency matrix $A$ for a graph and the ${\rm level} = n-1$ as input 
arguments, returns the path length matrix associated with the adjacency matrix $A$.
The function implements the dynamic programming Bellmann-Ford algorithm for solving 
the well-known ``all-pairs shortest path problem''. The algorithm requires 
${\cal O}(n^2m)$ arithmetic floating point operations (flops), where $n$ is the 
number of vertices and $m$ is the number of edges of the graph; if the graph is undirected, 
then the cost of the algorithm is halved. Notice that the function can be applied to 
determine shortest paths in a weighted graph having positive or negative weights.

\begin{algorithm}[ht]
\begin{algorithmic}[1]
\Function{\text{pathlength\_matrix}}{$A,\text{level}$}
\State $n=\text{size}(A,1)$;
\For {$i = 1:n$}
  \For {$j = 1:n$}
     \If{$A(i,j) = 0$ and $i \ne j$}
        \State $A(i,j) = \text{inf}$;
     \EndIf  
  \EndFor 
\EndFor 
\State $B = A$;
\State $C = \text{inf}(n)$;
\For{$\text{count} = 1:\text{level}-1$}
   \For{$i = 1:n$}
        \For{$j = 1:n$}
           \If{$j\ne i$}
               \For{$k = 1:n$}
                  \If{$B(k,j) \ne \text{inf}$ and $k \ne j$}
                     \State $C(i,j) = \min(C(i,j), A(i,k) + B(k,j))$; 
                  \EndIf
               \EndFor
            \EndIf 
        \EndFor
    \EndFor
    \State $A = C$; 
    \For{$i = 1:n$}
        \State $A(i,i) = 0$;
    \EndFor
\EndFor\\
\Return \text{Path length matrix};
\EndFunction
\end{algorithmic}
\end{algorithm}

In line 2 of the function Pathlength\_matrix, the function call ${\sf size}(A,1)$ 
yields the order of the matrix $A\in\R^{n\times n}$; in line
6, the matrix entry $A(i,j)$ is set to $\infty$. Similarly, in line 11, $C$ is defined as 
an $n\times n$ matrix with all entries equal to $\infty$. If the graph is unweighted and 
connected, then the above MATLAB function can be modified by introducing a \rm{break} 
before the last ${\bf end~for}$ when there is no entry $\infty$. The diameter of the graph
then is \rm{count+1}.

For both weighted and unweighted graphs, also when the graph is not connected, the above 
MATLAB function with argument ${\rm level}=K$, where $1<K<n$, computes  $A^{K,\star}$, 
i.e., the matrix of the distances between any distinct vertices of the graph using paths 
with at most $K$ edges. We note that the triangle inequality might not hold for the 
entries of $A^{K,\star}$ and some entries of this matrix may have the value $\infty$.

We will see how the matrix $A^{n-1,\star}$ associated with a connected graph $\mathcal{G}$
sheds light on the communication within the network determined by the graph. In fact, as a 
measure of the ease of communication in the graph, we like to use the average inverse 
geodesic length of $\mathcal{G}$  (i.e., its {\it global efficiency}, cf.  Section 
\ref{s3}) instead of the maximum geodesic length of $\mathcal{G}$ (i.e., its diameter). To this 
end, we introduce the {\it reciprocal path length matrix} 
$A^{n-1,\star,-1}=[a^{(n-1,\star,-1)}_{ij}]_{i,j=1}^n$ obtained by replacing the 
off-diagonal entries of the path length matrix by their reciprocals, i.e., 
\[
a^{(n-1,\star,-1)}_{ij}=1/a^{(n-1,\star)}_{ij},\qquad 1\leq i,j\leq n,
\]
where $1/\infty$ is identified with $0$.

We are interested in determining the shortest paths that use at most $K$ edges. We 
therefore also consider the {\it reciprocal $K$-path length matrix} 
$A^{K,\star,-1}=[a^{(K,\star,-1)}_{ij}]_{i,j=1}^n$, with 
$a^{(K,\star,-1)}_{ij}=1/a^{(K,\star)}_{ij}$. Thus, the entry $a^{(K,\star,-1)}_{ij}$ 
vanishes if $a^{(K,\star)}_{ij}=\infty$. Note that the same would happen if $K=n-1$, in 
case the graph $\mathcal{G}$ that determines the adjacency matrix $A$ is not connected. 
The matrix $A^{K,\star,-1}$ allows us to define the {\it global $K$-efficiency} of 
$\mathcal{G}$; see Section \ref{s3}.

In order to enhance communication using paths with at most $K$ edges, with $1<K<n$, i.e., 
to increase the global $K$-efficiency of the graph associated with the adjacency matrix 
$A$, we select edge weights by analyzing centrality properties of the vertices of the 
graph and, if computationally feasible, the spectral properties of the reciprocal $K$-path
length matrix $A^{K,\star,-1}$. In detail, if $K\ll n$ and the (sparse) non-negative 
reciprocal $K$-path length matrix $A^{K,\star,-1}$ is irreducible, then we apply the 
Perron-Frobenius theory by following the approach in \cite{djnr_new,sm1}. In our context,
the choice of weights is dictated by the analysis of the sensitivity to perturbations in 
$A$.

Applications of our approach include city planning and information transmission. As for 
disease propagation, a recent research study of the Zhejiang City Planning Center (China) 
pointed out a strong connection between the spread of the Covid-19 epidemic and the shape 
of the city. Cities with a radial structure (such as Milan) have good internal connections 
due to the capillarity of the public transport system - buses, trams, subways, and trains.
This dynamic made citizens of such cities more vulnerable to the arrival of the Covid-19 
virus: the incidence of infections compared to the number of inhabitants was generally 
larger in cities with a radial urban structure than in cities without this structure, due to the
good communication of the people in cities with radial urban structure.

This paper is organized as follows: Section \ref{s2} analyzes differences and similarities
of powers and tropical powers of the adjacency matrix for undirected and unweighted 
graphs. Section \ref{s3} reviews well-known measures that can be easily computed by means 
of the path length matrix and introduces novel ones. In Section \ref{s4} we present two 
algorithms that determine which edge-weight should be changed in order to boost global 
efficiency. Changing the edge weights may entail widening streets or increasing the
number of trams on a route, decreasing travel times on a highway by increasing 
the travel speed, or decreasing the waiting time for trams on a route. Finally, 
numerical tests are reported in Section \ref{s7} and concluding remarks can be found in 
Section \ref{s8}.

\section{Powers versus tropical powers}\label{s2}
Consider an undirected and unweighted simple graph $\mathcal{G}$ with adjacency matrix 
$A\in\R^{n\times n}$. Then the entry $a^{(k)}_{ij}$ of the matrix 
$A^k=[a^{(k)}_{ij}]_{i,j=1}^n\in\R^{n\times n}$ counts the number of walks of length $k$
between the vertices $v_i$ and $v_j$. A \emph{matrix function} based on the powers $A^k$ 
that is analytic at the origin, and vanishes there, can be defined by a formal Maclaurin 
series
\begin{equation}\label{matfun}
f(A)=\sum_{k=1}^{\infty}c_kA^k,
\end{equation}
where we for the moment ignore the convergence properties of this series. Usually long 
walks are considered less important than short walks, because information flows more 
easily through short walks than through long ones. Therefore matrix functions applied in 
network analysis generally have the property that $0\leq c_{k+1}\leq c_k$ for all 
$k\geq1$. The most common matrix function used in network analysis is the matrix 
exponential; see \cite{DMR,DMR2,estrada2011structure,estrada2005subgraph,FRR} for 
discussions and illustrations. We prefer to use the the modified matrix exponential 
\[
\exp_0(A):=\exp(A)-I,
\]
where $I$ denotes the identity matrix, because the first term in the Maclaurin series of 
$\exp(A)$ has no natural interpretation in the context of network modeling. For the 
modified matrix exponential, we have $c_k=1/k!$, and the series \eqref{matfun} converges 
for any adjacency matrix $A$. 

The {\it communicability} between distinct vertices $v_i$ and $v_j$, $i\ne j$, is defined
by
\[
[\exp_0(A)]_{ij}=\sum_{k=1}^{\infty}\frac{a^{(k)}_{ij}}{k!};
\]
see \cite{estrada2011structure} for the analogous definition based on $\exp(A)$. The
communicability accounts for all possible routes of communication between the vertices 
$v_i$ and $v_j$ in the network defined by the adjacency matrix $A$, and assigns a larger 
weight to shorter walks than to longer ones. The larger the value of $[\exp_0(A)]_{ij}$, 
the better is the communicability between the vertices $v_i$ and $v_j$. 

\begin{remark}\label{nnz}
Notice that even if there exists an integer $\widehat{k}$, $1\leq\widehat{k}< n-1$, such 
that $a_{ij}^{(\widehat{k})}>0$, one may have $a_{ij}^{(\widehat{k}+1)}=0$. Information 
provided by $A^{n-1}$ does not include information provided by all $A^{h}$ for 
$1\leq h<n-1$. This is one of the reasons for the interest in the matrix functions
$\exp_0(A)$ and $\exp(A)$.
\end{remark}

It is straightforward to show the following result.

\begin{proposition}\label{fnnz}
Let $\widehat{k}$, $1\leq\widehat{k}< n$, be the smallest integer power such that 
$a_{ij}^{(\widehat{k})}=p>0$, that is to say, 
\[
a^{(h)}_{ij}=0,\quad \forall h<\widehat{k},\quad a^{(\widehat{k})}_{ij}=p.
\] 
Then $\widehat{k}$ is the length of the shortest path that connects $v_i$ and $v_j$, i.e.,
\[
a_{ij}^{(h,\star)}=\infty,\quad \forall h<\widehat{k},\quad a_{ij}^{(h,\star)} = 
\widehat{k},\quad \forall h\geq\widehat{k},
\] 
and $p$ is the number of shortest paths that connect $v_i$ and $v_j$. 
\end{proposition}

\begin{example}\label{ex1}
In view of Proposition \ref{fnnz}, the information provided by both the path length matrix
and a suitable power of the adjacency matrix is of interest. Consider the undirected and 
unweighted graphs $\mathcal{G}_1$ and $\mathcal{G}_2$ depicted in Figures 
\ref{fig1}-\ref{fig2}. The adjacency matrices of these graphs are 
\begin{equation*}
A_1=\begin{bmatrix} 
0 & 0 & 1 & 1 & 1 \\ 
0 & 0 & 1 & 1 & 1\\ 
1 & 1 & 0 & 0 & 0\\ 
1 & 1 & 0 & 0 & 0\\ 
1 & 1 & 0 & 0 & 0
\end{bmatrix},
\quad\quad
A_2=\begin{bmatrix} 
0 & 0 & 1 \\ 
0 & 0 & 1   \\ 
1 & 1 & 0 
\end{bmatrix}.
\end{equation*}

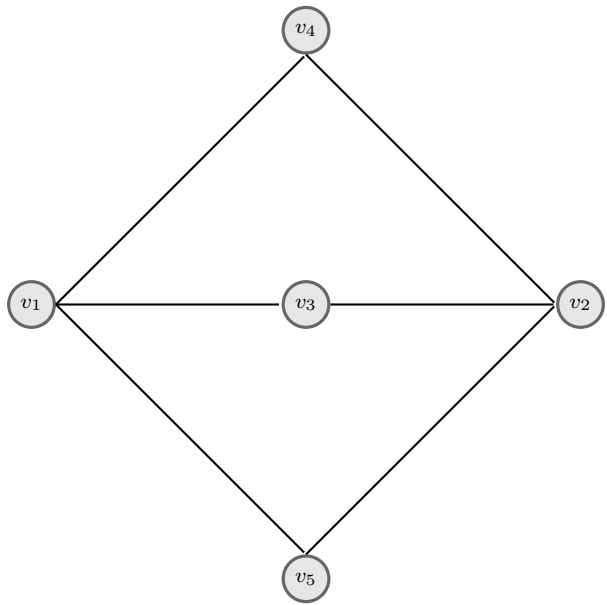
\begin{figure}[H]
\begin{center}
\begin{tikzpicture}[
 roundnode/.style={circle, draw=black!60, fill=black!10, very thick, 
 minimum size=6mm},
 ->,>=stealth',shorten >=1pt,node distance=2.5cm,auto,thick]
\node[roundnode]      (node1)                                 {$v_1$};
\node[roundnode]      (node3)       [right=3cm of node1]      {$v_3$};
\node[roundnode]      (node2)       [right=3cm of node3]      {$v_2$};
\node[roundnode]      (node4)       [above=3cm of node3]   {$v_4$};
\node[roundnode]      (node5)       [below=3cm of node3]    {$v_5$};
\draw[-] (node1.0) -- (node3.180);
\draw[-] (node3.0) -- (node2.180);
\draw[-] (node1.0) -- (node4.-90);
\draw[-] (node1.0) -- (node5.90);
\draw[-] (node4.-90) -- (node2.180);
\draw[-] (node5.90) -- (node2.180);
\end{tikzpicture}
\caption{Graph $\mathcal{G}_1$ in Example \ref{ex1}.}
\label{fig1}
\end{center}
\end{figure}

\begin{figure}[H]
\begin{center}
\begin{tikzpicture}[
 roundnode/.style={circle, draw=black!60, fill=black!10, very thick, minimum
 size=6mm},
 ->,>=stealth',shorten >=1pt,node distance=2.5cm,auto,thick]
\node[roundnode]      (node1)                                 {$v_1$};
\node[roundnode]      (node3)       [right=3cm of node1]      {$v_3$};
\node[roundnode]      (node2)       [right=3cm of node3]      {$v_2$};
\draw[-] (node1.0) -- (node3.180);
\draw[-] (node3.0) -- (node2.180);
\end{tikzpicture}
\caption{Graph $\mathcal{G}_2$ in Example \ref{ex1}.} 
\label{fig2}
\end{center}
\end{figure}

For both graphs, the shortest path between the vertices $v_1$ and $v_2$ has length $2$. 
However, in $\mathcal{G}_1$ there are three shortest paths that connect $v_1$ and $v_2$, 
while there is only one in $\mathcal{G}_2$. Thus, there is surely better communication 
between these vertices in $\mathcal{G}_1$ than in $\mathcal{G}_2$, even though this 
information is not provided by the path length matrices for these graphs. The path length
matrices are 
\begin{equation}\label{plm}
A_1^{4,\star}=\begin{bmatrix} 
0 & 2 & 1 & 1 & 1 \\ 
2 & 0 & 1 & 1 & 1\\ 
1 & 1 & 0 & 2 & 2\\ 
1 & 1 & 2 & 0 & 2\\ 
1 & 1 & 2& 2 & 0
\end{bmatrix},
\quad\quad
A_2^{2,\star}=\begin{bmatrix} 
0 & 2 & 1 \\ 
2 & 0 & 1   \\ 
1 & 1 & 0  
\end{bmatrix}.
\end{equation}
On the other hand, the second powers of the above adjacency matrices are 
\begin{equation*}
A_1^2=\begin{bmatrix} 
3 & 3 & 0 & 0 & 0 \\ 
3 & 3 & 0 & 0 & 0\\ 
0 & 0 & 2 & 2 & 2\\ 
0 & 0 & 2 & 2 & 2\\
0 & 0 & 2 & 2 & 2
\end{bmatrix},
\quad\quad
A_2^2=\begin{bmatrix} 
1 & 1 & 0 \\ 
1 & 1 & 0 \\ 
0 & 0 & 2 \\ 
\end{bmatrix}.
\end{equation*}
\end{example}

\begin{remark} For weighted graphs, the interpretation of the entry $a^{(k)}_{ij}$ of the 
matrix $A^k$ has to be modified. Indeed, $a^{(k)}_{ij}$ yields the sum of all products of 
all weights of the edges in the walks of length $k$ between the vertices $v_i$ and $v_j$. 
(In the unweighted case, any product of all weights of the edges in a walk is $1$.) Remark 
\ref{nnz} holds true, while Proposition \ref{fnnz} does not.
\end{remark}

\section{Measures that depend on the path length matrix}\label{s3}
Let for now the graph $\mathcal{G}$ be undirected.

\subsection{Closeness centrality} 
Let the graph $\mathcal{G}$ be connected. Then the reciprocal of the sum of all lengths of 
the shortest paths starting from vertex $v_i$, 
\[
c_i=\frac{1}{\sum_{j\ne i}a_{ij}^{(n-1,\star)}}=[A^{n-1,\star}\mathbf{1}_i]^{-1},
\]
where $\mathbf{1}_i\in\R^n$ denotes the vector with all zero entries except for the 
$i^{th}$ entry, which is one, is referred to as the {\it closeness centrality} of $v_i$ in
$\mathcal{G}$; see, e.g., \cite{BBV}. This measure gives a large value to vertices that 
have small shortest path distances to the other vertices of the graph.

\subsection{The radius and center of a network}  
The maximum length over all shortest paths starting from vertex $v_i$ in a weighted 
or unweighted connected graph ${\mathcal{G}}$, that is to say,
\[
e_i=\max_{j\ne i} a_{ij}^{(n-1,\star)},
\]
is commonly referred to as the {\it eccentricity} of the vertex $v_i$. The diameter may be seen 
as the maximum eccentricity among the vertices of the network. The radius $r_{\mathcal{G}}$ 
of $\mathcal{G}$ is the minimum eccentricity of a vertex. One has 
\[
r_{\mathcal{G}}=\min_i\max_{j\ne i} a_{ij}^{(n-1,\star)}.
\] 
It is shown in \cite{CH} that
\[
r_{\mathcal{G}}\leq d_{\mathcal{G}}\leq 2r_{\mathcal{G}}.
\] 
A vertex is said to be \emph{central} if its eccentricity is equal to the radius of the 
graph. The \emph{center} of the graph is the set of all central vertices; see, e.g., 
\cite{estrada2011structure}.

\subsection{Average shortest path length}  
The average shortest path length of a weighted or unweighted connected graph $\mathcal{G}$ 
computed over all possible pairs of vertices in the network \cite{BBV} is given by 
\[
a_{\mathcal{G}}=\frac{1}{n(n-1)}{\sum_{i,j\ne i} a_{ij}^{(n-1,\star)}}=
\frac{1}{n(n-1)}\mathbf{1}^TA^{n-1,\star} \mathbf{1},
\] 
where $\mathbf{1}\in\R^n$ denotes the vector with all entries one. Let the graph 
$\mathcal{G}$ be unweighted and be formed by a path of $n$ vertices. Then the largest 
average shortest path length is $a_{\mathcal{G}}=(n+1)/3$. If $a_{\mathcal{G}}$ scales 
logarithmically with $n$, then $\mathcal{G}$ displays the ``small-world phenomenon''; 
see, e.g., \cite{estrada2011structure}.

\subsection{Harmonic centrality and global efficiency} 
The \emph{efficiency} of a path between any two vertices of a weighted or unweighted graph 
$\mathcal{G}$ is defined as the inverse of the length of the path. The sum of the inverses
of the length of all shortest paths starting from vertex $v_i$, i.e., the sum of the 
efficiencies of all shortest paths starting from $v_i$,
\begin{equation}\label{h_i}
h_i=\sum_{j\ne i}[a_{ij}^{(n-1,\star)}]^{-1},
\end{equation}
is referred to as the \emph{harmonic centrality} of $v_i$; see, e.g., \cite{BBV}. The 
latter measure gives a large centrality to vertices $v_i$ that have small shortest path 
distances to the other vertices of the graph. Harmonic centrality will play a central role in
our analysis. We define the \emph{$h$-center} of a graph as the set of all vertices with 
the largest harmonic centrality.

If the graph is connected, then the average shortest path efficiency over all possible 
pairs of vertices, also known as the \emph{average inverse geodesic length}, is referred 
to as the \emph{global efficiency} of the graph \cite{BBV}:
\begin{equation}\label{e_G}
e_{\mathcal{G}}=\frac{1}{n(n-1)}{\sum_{i,j\ne i} [a_{ij}^{(n-1,\star)}]^{-1}}.
\end{equation}

Finally, we remark that in the context of molecular chemistry, the sum of reciprocals of 
distances between all pairs of vertices of an undirected and unweighted connected graph is
known as the \emph{Harary index} and the reciprocal path length matrix $A^{n-1,\star,-1}$ 
is referred to as the \emph{Harary matrix} \cite{PNTM}. 

\begin{remark}
The measures \eqref{h_i} and \eqref{e_G} can be especially useful when the network has
more than one connected component, because infinite distances do not contribute to these
``harmonic'' averages.
\end{remark}

\begin{example}\label{ex2}
As an illustration of the above measures, consider again the graphs $\mathcal{G}_1$ and 
$\mathcal{G}_2$ of Example \ref{ex1}. Table \ref{tab1} reports the diameter, radius, 
average shortest path length, and global efficiency of these graphs. Notice that the 
center of $\mathcal{G}_1$ is given by the set of all vertices and the center of 
$\mathcal{G}_2$ is made up of the vertex $v_3$, only. The $h$-center of $\mathcal{G}_1$ is
formed by the vertices $v_1$ and $v_2$, and the $h$-center of $\mathcal{G}_2$ is given by 
the vertex $v_3$, only. Table \ref{tab2} shows the eccentricity, harmonic centrality, and
closeness centrality of all the vertices of $\mathcal{G}_1$ and 
$\mathcal{G}_2$. All the measures in this example are computed by using the path 
length matrices \eqref{plm}.
\end{example}

\begin{table}
\centering
\begin{tabular}{l*{4}{|c}}
$\mathcal{G}$           & $d_\mathcal{G}$ & $r_\mathcal{G}$ & 
$a_\mathcal{G}$& $e_\mathcal{G}$ \\
\hline
$\mathcal{G}_1$            & 2 & 2 & 1.40 &  0.80\\
$\mathcal{G}_2$           & 2 & 1 & 1.33 &  0.83
\end{tabular}
\caption{Example \ref{ex2}. Diameter, radius, average shortest path length, and global 
efficiency of the graphs $\mathcal{G}_1$ and $\mathcal{G}_2$ depicted in Figures 
\ref{fig1}-\ref{fig2}. Note that the average shortest path length of $\mathcal{G}_2$ is
maximal; indeed, $(n+1)/3=1.33$.} 
\label{tab1}
\end{table}

\begin{table}[ht]
\centering
\begin{tabular}{c*{3}{|c}}
$v_i (\mathcal{G}_1)$           & $e_i$ & $h_i$ & $c_i$  \\
\hline
1           & 2 & 3.50 & 0.20   \\
2           & 2 & 3.50 & 0.20\\
3           & 2 & 3.00 &  0.17\\
4           & 2 & 3.00 &  0.17\\
5           & 2 & 3.00 &  0.17
\end{tabular}
\qquad\qquad
\begin{tabular}{c*{3}{|c}}
$v_i(\mathcal{G}_2)$           & $e_i$ & $h_i$ & $c_i$  \\
\hline
1          & 2 & 1.50 & 0.33   \\
2          & 2 & 1.50 &  0.33\\
3         & 1 & 2.00 &  0.50
\end{tabular}

\caption{Example \ref{ex2}. Eccentricity, harmonic centrality, and closeness centrality of 
the vertices of  $\mathcal{G}_1$ (left table) and of  $\mathcal{G}_2$ (right table).} 
\label{tab2}
\end{table}

\subsection{Harmonic $K$-centrality and global $K$-efficiency}
When one considers shortest paths that are made up of at most $K$ edges, the matrix 
$A^{K,\star,-1}$ takes the role of the reciprocal path length matrix $A^{n-1,\star,-1}$.
We are in a position to introduce the \emph{harmonic $K$-centrality} of the vertex $v_i$. 
It is given by
\begin{equation*}\label{h_i_K}
h_i^{K}=\sum_{j\ne i}[a_{ij}^{(K,\star)}]^{-1}.
\end{equation*}
The {\it global $K$-efficiency} of a graph $\mathcal{G}$ is defined by
$$
e_{\mathcal{G}}^{K}=\frac{1}{n(n-1)}{\sum_{i,j\ne i} [a_{ij}^{(K,\star)}]^{-1}}.
$$
The set of all vertices of a graph $\mathcal{G}$ with the largest harmonic $K$-centrality
is referred to as the {\it $h^{K}$-center} of the graph.

\subsection{Out-centrality versus in-centrality}
Let the  graph $\mathcal{G}$ be directed. Then the above centrality measures (closeness,
eccentricity, harmonic centrality, and harmonic $K$-centrality) determine the importance 
of a vertex $v_j$ by taking into account the paths that start at $v_j$. 
These measures therefore may be considered measures of out-centrality. One also may be 
interested in measuring the importance of a vertex $v_j$ by considering the paths that end at it, that 
is to say by measuring the in-centrality of $v_j$. This can be achieved by replacing the 
path length matrix $A^{n-1,\star}$ in the measures mentioned by its transpose. This allows
us to introduce the measures \emph{closeness in-centrality}, \emph{in-eccentricity}, 
\emph{harmonic in-centrality}, and \emph{harmonic $K_{\rm in}$-centrality}. These measures
are defined as 
$$
c_j^{\rm in}=\frac{1}{\sum_{i\ne j}a_{ij}^{(n-1,\star)}} \,,\;\;\; e_j^{\rm in}=
\max_{i\ne j} a_{ij}^{(n-1,\star)}\,, \;\;\; h_j^{\rm in}=
\sum_{i\ne j}[a_{ij}^{(n-1,\star)}]^{-1}\,,\;\;\; h_j^{K_{\rm in}}=
\sum_{i\ne j}[a_{ij}^{(K,\star)}]^{-1}.
$$
We also define the {\it in-radius} of  $\mathcal{G}$,
$$
r^{\rm in}_{\mathcal{G}}=\min_j\max_{i\ne j} a_{ij}^{(n-1,\star)}.
$$
Finally, the vertex $v_j$ is said to be \emph{in-central} if its in-eccentricity equals 
the in-radius of the graph. The \emph{$h^{K_{\rm in}}$-center} of a graph $\mathcal{G}$ is
the set of all vertices with the largest harmonic $K_{\rm in}$-centrality. The notions of 
$h^{K_{\rm in}}$-center and $h^{K_{\rm out}}$-center will be of interest in the sequel.

\section{Enhancing network communication} \label{s4}
The diameter of a weighted or unweighted graph provides a measure of how easy it is for 
the vertices of the graph to communicate. When graphs are used as models for communication
networks, the diameter plays an important role in the performance analysis and cost 
optimization. A simple way to decrease the diameter of a graph so that information can be 
transmitted more easily between vertices of the graph \cite{djnr_new} is to decrease the 
weight of an edge that belongs to all maximal shortest paths, if feasible. 

\begin{example}\label{ex3}
Consider the graph $\mathcal{G}_2$ in Example \ref{ex1}. The edge $e(v_2\leftrightarrow v_3)$ belongs to the
maximal shortest path and is a {\it bridge}, that is its removal would make the
vertex $v_2$ unreachable and the perturbed graph $\tilde{\mathcal{G}}_2$ disconnected. In 
detail, the edge $e(v_2\leftrightarrow v_3)$ belongs to the shortest path between the
vertices $v_1$ and $v_2$, which is the only 
shortest path of length $d_\mathcal{G}$. Thus, decreasing the weight of the edge 
$e(v_2\leftrightarrow v_3)$ decreases both the average shortest path length and the 
diameter of the graph. Specifically, if one decreases the weights $a_{23}$ and $a_{32}$ 
from $1$ to $0.5$, one obtains the matrices
\begin{equation*}
\tilde{A}_2=\begin{bmatrix} 
0 & 0 & 1 \\ 
0 & 0 & 0.5   \\ 
1 & 0.5 & 0 
\end{bmatrix}, 
\quad\quad
\tilde{A}_2^{2,\star}=\begin{bmatrix} 
0 & 1.5 & 1 \\ 
1.5& 0 & 0.5  \\ 
1 & 0.5 & 0 
\end{bmatrix}.
\end{equation*}
The perturbed graph $\tilde{\mathcal{G}}_2$ has diameter $d_{{\mathcal{\tilde G}}_2}=1.5$, 
radius $r_{{\mathcal{\tilde G}}_2}=1$, average shortest path length 
$a_{{\mathcal{\tilde G}}_2}=1$, and global efficiency $e_{{\mathcal{\tilde G}}_2}=1.22$. 
The eccentricity, harmonic centrality, and closeness centrality of the vertices of the 
graph $\tilde{\mathcal{G}}_2$ are reported in Table \ref{tab3}.

\begin{table}
\centering
\begin{tabular}{c*{3}{|c}}
$v_i(\mathcal{\tilde{G}}_2)$           & $e_i$ & $h_i$ & $c_i$  \\
\hline
1          & 1.50 & 1.67 & 0.40   \\
2          & 1.50 & 2.67 &  0.50\\
3         & 1.00 & 3.00 &  0.67
\end{tabular}

\caption{Example \ref{ex3}. Eccentricity, harmonic centrality, and closeness 
centrality of the vertices of  $\mathcal{\tilde{G}}_2$.}\label{tab3}
\end{table}
\end{example}

\subsection{Increasing the global $K$-efficiency}
We propose two approaches to increase the global efficiency of a network.  

\subsubsection{The function eKG1}
Let for now the graph $\mathcal{G}$ be directed. The first approach is based on the 
observation that the most important vertices with respect to the global $K$-efficiency 
live in the vertex subsets $h^{K_{\rm out}}$-center and $h^{K_{\rm in}}$-center of the 
graph. These
vertices may be interpreted as important intermediaries, that quickly collect information 
from many vertices and quickly broadcast it to many others vertices. Indeed, strengthening 
an existing connection from a vertex of the $h^{K_{\rm in}}$-center to a vertex of the 
$h^{K_{\rm out}}$-center is likely to  strengthen their communicability by having new 
shorter paths with at most $K$ steps that exploit these connections. This is likely to 
increase the global $K$-efficiency more than strengthening an existing connection between 
vertices with lower harmonic $K_{\rm in}$- and $K_{\rm out}$-centrality. Here we consider 
graphs whose edge weights represent travel times or waiting times. Hence, strengthening is
achieved by \emph{decreasing} appropriate weights. 

We construct the perturbed adjacency matrix
\begin{equation}\label{Atilde1}
 \tilde{A}=A+\gamma_{h_1,h_2}\mathbf{1}_{h_1}\mathbf{1}_{h_2}^T, 
 \;\;\; \mbox{with}\;\;\gamma_{h_1,h_2} = - a_{h_1,h_2}/2,
\end{equation}
 
If the graph is undirected, the above approach simplifies, because the $h^{K_{\rm out}}$-
and $h^{K_{\rm in}}$-centers coincide. Hence, the idea is to strengthen the connection 
between vertices with the largest harmonic $K$-centrality. The perturbed adjacency matrix 
$\tilde{A}$ will be
\begin{equation}\label{Atilde2}
 \tilde{A}=A+\gamma_{h_1,h_2}(\mathbf{1}_{h_1}\mathbf{1}_{h_2}^T +
 \mathbf{1}_{h_2}\mathbf{1}_{h_1}^T) .
\end{equation}
The MATLAB function eKG1 describes the necessary computations. 
The operator $==$ in line 5 of the function eKG1 stands for logical equal to, and the 
symbol $./$ in line 6 denotes element-wise division. The function call ${\sf sum}(M)$ for 
a matrix $M\in\R^{n\times n}$ computes a row vector $m\in\R^n$, whose $j^{th}$ component 
is the 
sum of the entries of the $j^{th}$ column; the function call ${\sf sum}(M,2)$ for a 
matrix $M\in\R^{n\times n}$ computes a column vector, whose $i^{th}$ entry is the sum
of the elements of row $i$ of $M$. The blip in line 8 denotes transposition, and the 
operator $.*$ in line 11 stands for vector-vector element-wise product.

\begin{algorithm}[ht]
\begin{algorithmic}[1]
\Function{\text{eKG1}}{$A,\text{flagsym},K$};
\State $n=\text{size}(A,1)$;
\State $P=\text{pathlength\_matrix}(A,K)$;        
\State $Pr=P$;
\State $Pr(Pr==\text{inf})=0$;
\State $Pr(Pr>0)=1./Pr(Pr>0)$;            
\State $H_{\text{in}}=\text{sum}(Pr)$;
\State $H=\text{sum}(Pr,2)'$;
\State $\text{eKG}=1/(n(n-1))\,\text{sum}(H)$;          
\State $[\sim,h]=\max(H_{in})$;
\State $[\sim,k]=\max(H.*A(h,:))$;            
\State $A(h,k)=A(h,k)/2$; 
\If{\text{flagsym}}
\State $A(k,h)=A(h,k)$;
\EndIf\\
\Return $A$, $\text{eKG}$; 
\EndFunction 
\end{algorithmic}
\end{algorithm}

\subsubsection{The function eKG2}
Assume for now that the graph $\mathcal{G}$ is directed. If computing the Perron root 
$\rho_K$ and the unique positive left and right eigenvectors of unit norm (the Perron 
vectors) of the reciprocal $K$-path length matrix $A^{K,\star,-1}$ is not computationally 
feasible or if this matrix is not irreducible, then one can use the function eKG1. 
However, if $A^{K,\star,-1}$ is irreducible and its left and right Perron vectors 
$\mathbf{x}_K=(x_{K,i})$ and $\mathbf{y}_K=(y_{K,i})$ can be computed, then these vectors
determine the Wilkinson perturbation $W_K=\mathbf{y}_K\mathbf{x}_K^T$; see 
\cite[Section 2]{W}. Following \cite{sm1}, to induce the maximal perturbation in 
$\rho_K$, one chooses the indices $(h_1,h_2)$ such that $W_K(h_1,h_2)$ is the largest 
entry of $W_K$ and $A(h_1,h_2)>0$, i.e., the indices of the largest entry of the Wilkinson
perturbation ``projected'' onto the zero-structure of $A$; see, e.g. \cite{NP06}. Thus,
\[
(h_1,h_2):  x_{K,h_2}\,y_{K,h_1}=(\mathbf{y}_K\mathbf{x}_K^T)_{h_1, h_2}=
\max_{i,j:A(i,j)>0} (W_K)_{i,j}.
 \]
As in function eKG1, one strengthens the edge $e(v_{h_1}\rightarrow v_{h_2})$ by halving
its weight. The perturbed adjacency matrix then is given by \eqref{Atilde1}.

When the graph $\mathcal{G}$ is undirected, the left and right Perron vectors coincide and
the perturbed adjacency matrix is constructed as in \eqref{Atilde2}. The outlined approach
is implemented by the MATLAB function eKG2. We recall that the function ${\sf abs}(\sf v)$ of 
a vector ${\sf v}$, used in the function eKG2, returns a vector, whose components are the 
absolute value of the components of $\sf v$. The expression ${\sf sum}({\sf sum}(Pr))$ on
line 7 sums all entries of the matrix $Pr$. The operator $.*$ on line 20 denotes the
Hadamard product of two matrices; the entries of the matrix $A>0$ are one if the 
corresponding entry of $A$ is positive; they are zero otherwise. The function 
{\sf ind2sub} determines the equivalent subscript values corresponding to a given single 
index into an array.  

\begin{algorithm}[ht]
\begin{algorithmic}[1]
\Function{\text{eKG2}}{$A,\text{flagsym},K$}
\State $n=\text{size}(A,1)$;
\State $P=\text{pathlength\_matrix}(A,K)$;        
\State $Pr=P$;
\State $Pr(Pr==\text{inf})=0$;
\State $Pr(Pr>0)=1./Pr(Pr>0)$;            
\State $\text{eKG}=1/(n(n-1))\,\text{sum}(\text{sum}(Pr))$;    
\If{\text{not~flagsym}}
   \State $[X,D,Y] = \text{eig}(Pr)$;
   \State $[\rho,\text{ind}]=\max(\text{abs}(\text{diag}(D)))$;
   \State $x=X(:,\text{ind})$;
   \State $y=Y(:,\text{ind})$
\Else
   \State $[X,D] = \text{eig}(Pr)$; 
   \State $[\rho,\text{ind}]=\max(\text{abs}(\text{diag}(D)))$; 
   \State $x=X(:,\text{ind})$;             
   \State $y=x$;
\EndIf
\State $W=yx'$; 
\State $W=W.*(A>0)$;
\State $[\sim,\text{lin}]=\max(W(:))$;
\State $[h,k]=\text{ind2sub}([n,n],\text{lin})$;         
\State $A(h,k)=A(h,k)/2$; 
\If{\text{flagsym}}
\State $A(k,h)=A(h,k)$;
\EndIf\\
\Return $A$, $\text{eKG}$
\EndFunction
\end{algorithmic}
\end{algorithm}

We expect the global $K$-efficiency to increase the most when decreasing the edge-weight 
that makes the Perron root $\rho_K$ change the most. The perturbation of the Perron root
generated by the function eKG2 typically is larger than the perturbation determined by 
the function eKG1. The difference in these perturbations is analogous to the difference 
between considering the most important vertex in a graph the one with the largest degree 
and the one with maximal eigenvector centrality.

\begin{remark}
Both functions eKG1 and eKG2 maximize lower bounds for the global $K$-efficiency of the 
network. Consider for the sake of clarity the undirected case. Let $\mathbf{h}_K$ denote
the vector of the harmonic $K$-centralities of the vertices of the graph. Its $1$-norm is
the sum in the numerator of the global $K$-efficiency, and its $\infty$-norm is what 
function eKG1 is maximizing. Indeed, one has
\[
n(n-1)\,e_{\mathcal{G}}=\|\mathbf{h}_K\|_1\geq \|\mathbf{h}_K\|_{\infty}=
\|A^{K,\star,-1}\|_{\infty} \geq\rho_K.
\]
\end{remark}

\begin{remark} 
The Perron root is bounded from below and from above by the minimal and maximal entries of  
$\mathbf{h}_K$, respectively. We have
\begin{equation}\label{rho}
\min_{i}\sum_{j\ne i}[a_{ij}^{(K,\star)}]^{-1}\leq\rho_K\leq 
\max_{i}\sum_{j\ne i}[a_{ij}^{(K,\star)}]^{-1}=
\|\mathbf{h}_K\|_{\infty}.
\end{equation}
Indeed, one has 
\[
\mathbf{1}^T A^{K,\star,-1} \mathbf{x}_K = \mathbf{1}^T \rho_K \mathbf{x}_K = 
\rho_K \|\mathbf{x}_K\|_1
\]
and
\[
\mathbf{1}^T A^{K,\star,-1} \mathbf{x}_K =  \mathbf{x}_K^T A^{K,\star,-1} \mathbf{1} 
=\sum_{i=1}^n(x_{K,i} \sum_{j \ne i} [a_{ij}^{(K,\star)}]^{-1} ),
\]
so that one obtains \eqref{rho} by observing that
\[
\min_{i}\sum_{j\ne i}[a_{ij}^{(K,\star)}]^{-1}\|\mathbf{x}_K\|_1\leq
\rho_K\|\mathbf{x}_K\|_1\leq 
\max_{i}\sum_{j\ne i}[a_{ij}^{(K,\star)}]^{-1}\|\mathbf{x}_K\|_1.
\]
\end{remark} 
 
\begin{example}\label{ex4}
We apply the functions eKG1 and eKG2 to the graphs of Example \ref{ex1}. For neither graph, 
there is a unique ``best choice'' to report. 
 
First consider the graph $\mathcal{G}_1$ in Example \ref{ex1}. Then 
$A_1^{2,\star}=A_1^{3,\star}= A_1^{4,\star}$. We let $K=2$ and obtain
\begin{equation*}
A_1^{2,\star,-1}=\begin{bmatrix} 
0 & 0.5 & 1 & 1 & 1 \\ 
0.5 & 0 & 1 & 1 & 1\\ 
1 & 1 & 0 & 0.5& 0.5\\ 
1 & 1 & 0.5 & 0 & 0.5\\ 
1 & 1 & 0.5& 0.5 & 0
\end{bmatrix}\,.
\end{equation*}
The vector $\mathbf{h}_2$ of harmonic $2$-centralities is 
$\mathbf{h}_2=[3.5,3.5,3,3,3]^T$, while the Perron vector $\mathbf{x}_2$ is given by
$\mathbf{x}_2=[0.47,0.47,0.43,0.43,0.43]^T$. This tells us that the  vertices $v_1$ and 
$v_2$ are the most important ones in the sense of both harmonic centrality and eigenvector
centrality. Indeed, these vertices are the only ones that are connected by paths of 
minimal length with three vertices. Thus, they are well connected. However, $A_1(1,2)=0$. 

Both functions eKG1 and eKG2 give $(h_1,h_2)=(1,3)$, and yield the matrix
\begin{equation*}
\tilde{A}_1=\begin{bmatrix} 
0 & 0 & 0.5 & 1 & 1 \\ 
0 & 0 & 1 & 1 & 1\\ 
0.5 & 1 & 0 & 0 & 0\\ 
1 & 1 & 0 & 0 & 0\\ 
1 & 1 & 0 & 0 & 0
\end{bmatrix}\,,
\end{equation*}
even though other choices of existing edges are equally valid. The global $2$-efficiency 
$e_{\tilde{\mathcal{G}}_1}^2$ of $\tilde{A}_1$ is $0.95$; compare with the global 
$2$-efficiency $e_{\mathcal{G}_1}^2=0.80$ of $A_1$. We remark that a
perturbation of another existing edge would have led to the same increase of the global 
$2$-efficiency.

Consider the graph $\mathcal{G}_2$ in Example \ref{ex1}. One has
\begin{equation*}
A_2^{2,\star,-1}=\begin{bmatrix} 
0 & 0.5 & 1  \\ 
0.5 & 0 & 1 \\ 
1 & 1 & 0 \\ 
\end{bmatrix}\,.
\end{equation*}
The vector of harmonic $2$-centralities is $\mathbf{h}_2=[1.5,1.5,2]^T$, while the 
Perron vector is given by $\mathbf{x}_2=[0.54,  0.54, 0.64]^T$. Both functions eKG1 and 
eKG2 yield $(h_1,h_2)=(3,1)$, even though the choice $(h_1,h_2)=(3,2)$ is equally valid
and gives the matrix
\begin{equation*}
\tilde{A}_2=\begin{bmatrix} 
0 & 0 & 0.5  \\ 
0 & 0 & 1 \\ 
0.5 & 1 & 0 \\ 
\end{bmatrix}.
\end{equation*}
The global $2$-efficiency of $\tilde{A}_2$ is $e_{\tilde{\mathcal{G}}_2}^2=1.22$, while 
the global $2$-efficiency of the matrix $A_2$ is $e_{\mathcal{G}_2}^2=0.83$. The diameter 
associated with the graph determined by $\tilde{A}_2$ is only $1.5$, while the diameter of
the graph associated with $A_2$ is $2$. The choice 
$(h_1,h_2)=(3,2)$, which was considered in Example \ref{ex3}, gives the same results.
\end{example}

\section{Numerical tests} \label{s7} 
The numerical tests reported in this section have been carried out using MATLAB R2022b on 
a $3.2$ GHz Intel Core i7 6 core iMac. The Perron root and left and right Perron vectors 
for small to moderately sized graphs can easily be evaluated by using the MATLAB function 
{\sf eig}. For large-scale graphs these quantities can be computed by the MATLAB function 
{\sf eigs} or by the two-sided Arnoldi algorithm, introduced by Ruhe \cite{R} and improved
by Zwaan and Hochstenbach \cite{ZH}. 

\begin{example}\label{ex5}
Consider the adjacency matrix for the network {\em Air500} in \cite{air500}. This data set 
describes flight connections for the top 500 airports worldwide based on total passenger 
volume. The flight connections between airports are for the year from 1 July 2007 to 30 
June 2008. The network is represented by a directed unweighted connected graph with $500$ 
vertices and $24009$ directed edges. The vertices of the network are the airports and the 
edges represent direct flight routes between two airports. 

The path length matrix $A^{5,\star}=A^{499,\star}$ yields the diameter and the radius of 
the graph $5$ and $3$, respectively. The information provided by the vector 
of the harmonic centralities and the Perron vector for the reciprocal path length matrix 
$A^{5,\star,-1}$ is the same as the one given by the vector of harmonic $K$-centralities 
and the Perron vector for $A^{K,\star,-1}$ with $K=2$; cf. Table \ref{Air500}. Therefore, 
the perturbation that increases the global $K$-efficiency the most also will enhance the 
global efficiency the most. The information provided by Table \ref{Air500} suggests that
the number of flights from the Frankfurt FRA Airport (vertex $v_{161}$) to the JFK Airport
in New York (vertex $v_{224}$) should be doubled in order to half the wait time between 
these flights. Doubling the number of flights corresponds to halving the weight for
the corresponding edge.

\end{example}

\begin{table}[ht]
\centering
\begin{tabular}{c*{3}{|c}}
$K$           & $(h_1,h_2)$ & $e_{\mathcal{G}}^K$ & $e_{\tilde{\mathcal{G}}}^K$  \\
\hline
5          & $(161,224)$ & 0.4839 & 0.4856   \\
4          & $(161,224)$ & 0.4839 &  0.4855\\
3         & $(161,224)$ & 0.4791 &  0.4807\\
2          & $(161,224)$ & 0.3604 &  0.3606
\end{tabular}

\caption{Example \ref{ex5}. Indices chosen by the functions eKG1 and eKG2 and the global 
$K$-efficiency of both the given graph, ${\mathcal{G}}$, and the perturbed graph, 
${\tilde{\mathcal{G}}}$, for $K=2,3,4,5$.} 
\label{Air500}
\end{table}

\begin{example}\label{ex6}
This example considers an undirected unweighted connected graph ${\mathcal{G}}$ that 
represents the German highway system network {\em Autobahn}. The graph is available
at \cite{air500}. Its $1168$ vertices are German locations and its $1243$ edges represent 
highway segments that connect them. 

Let $A$ be the adjacency matrix associated with ${\mathcal{G}}$. The path length matrix 
$A^{62,\star}=A^{1167,\star}$ shows that the diameter and the radius of ${\mathcal{G}}$ 
are $62$ and $34$, respectively, whereas its global efficiency equals 
$6.7175\cdot 10^{-2}$. One notices that there is only one shortest path of length $62$,
which connects the vertices $v_{116}$ and $v_{1154}$. The diameter of the graph can be
decreased by halving the weight of the edges $a_{120,116}$ and $a_{116, 120}$ (since the
graph is undirected), because this is the unique edge that connects the vertex $v_{116}$ 
to the other vertices of the graph.

Let ${\hat{A}}$ denote the perturbed adjacency matrix and $\hat{\mathcal{G}}$ the 
corresponding graph. The global efficiency of $\hat{\mathcal{G}}$ is 
$6.7177\cdot 10^{-2}$ and the diameter is $61.5$. We turn to the application of 
the functions eKG1 and eKG2 to increasing the global efficiency. Table \ref{T_1} reports 
the global $K$-efficiency (for several values of $K$) of the graph ${\tilde{\mathcal{G}}}$
associated with the adjacency matrix $\tilde{A}$ obtained by halving both entries 
$a_{h_1,h_2}$ and $a_{h_2,h_2}$ of $A$, computed by the function eKG1. Table \ref{T_2} shows 
the global $K$-efficiency (for several values of $K$) of the graph ${\tilde{\mathcal{G}}}$ 
associated with the adjacency matrix $\tilde{A}$, computed by the function eKG2. Also in 
this example, the information provided by the vector of harmonic centralities and the 
Perron vector for the reciprocal path length matrix $A^{62,\star,-1}$ is exactly the same 
information that is provided by the vector of harmonic $K$-centralities and the Perron 
vector for $A^{K,\star,-1}$ with $K\geq4$. We note that the latter vectors are less
expensive to determine than the former. The information provided by both tables suggests
that one should double the width of the highway that connects the cities of Duisburg 
(vertex $v_{219}$) and Krefeld (vertex $v_{565}$) to half the travel time. These cities 
are 10 miles apart. Doubling the width corresponds to halving the weight associated
with the corresponding edge.

\begin{table}[ht]
\centering
\begin{tabular}{c*{3}{|c}}
$K$           & $(h_1,h_2)$ & $e_{\mathcal{G}}^K$ & $e_{\tilde{\mathcal{G}}}^K$  \\
\hline
62          & $(219,565)$ & $ 6.7175   \cdot 10^{-2}$ & $ 6.7559  \cdot 10^{-2}$    \\
52          & $(219,565)$ & $6.7166 \cdot 10^{-2}$ &  $6.7550 \cdot 10^{-2}$\\ 
42         & $(219,565)$ & $6.6965 \cdot 10^{-2}$ &  $6.7349 \cdot 10^{-2}$\\
32          & $(219,565)$ & $6.5105 \cdot 10^{-2}$ & $6.5485 \cdot 10^{-2}$   \\
22          & $(219,565)$ & $5.5674 \cdot 10^{-2}$ &  $5.6024 \cdot 10^{-2}$\\
12        & $(219 ,565)$ & $2.8426 \cdot 10^{-2}$ &  $2.8621 \cdot 10^{-2}$\\
\hline
5     & $(219 ,565)$ & $7.9991 \cdot 10^{-3}$ &  $ 8.0304\cdot 10^{-3}$\\
4     & $(219 ,565)$ & $6.1823 \cdot 10^{-3}$ &  $6.2019\cdot 10^{-3}$\\
3     & $(219,217)$ & $4.6017 \cdot 10^{-3}$ &  $4.6112\cdot 10^{-3}$\\
2         & $(219,217)$ & $3.2082 \cdot 10^{-3}$ &  $3.2124 \cdot 10^{-3}$\\
\end{tabular}
\caption{Example \ref{ex6}. Indices chosen by function eKG1 and global 
$K$-efficiency of both  ${\mathcal{G}}$ and ${\tilde{\mathcal{G}}}$, for 
$K=2,3,4,5$ and for $K=12:10:62$.} 
\label{T_1}
\end{table}

\begin{table}[ht]
\centering
\begin{tabular}{c*{3}{|c}}
$K$           & $(h_1,h_2)$ & $e_{\mathcal{G}}^K$ & $e_{\tilde{\mathcal{G}}}^K$  \\
\hline
62          & $(565 ,219)$ & $ 6.7175   \cdot 10^{-2}$ & $ 6.7559  \cdot 10^{-2}$    \\
52          & $(565 ,219)$ & $6.7166 \cdot 10^{-2}$ &  $6.7550 \cdot 10^{-2}$\\ 
42         & $(565 ,219)$ & $6.6965 \cdot 10^{-2}$ &  $6.7349 \cdot 10^{-2}$\\
32          & $(565 ,219)$ & $6.5105 \cdot 10^{-2}$ & $6.5485 \cdot 10^{-2}$   \\
22          & $(565 ,219)$ & $5.5674 \cdot 10^{-2}$ &  $5.6024 \cdot 10^{-2}$\\
12        & $(565 ,219)$ & $2.8426 \cdot 10^{-2}$ &  $2.8621 \cdot 10^{-2}$\\
\hline
5     & $(565 ,219)$ & $7.9991 \cdot 10^{-3}$ &  $ 8.0304\cdot 10^{-3}$\\
4     & $(565 ,219)$ & $6.1823 \cdot 10^{-3}$ &  $6.2019\cdot 10^{-3}$\\
3     & $(267,219)$ & $4.6017 \cdot 10^{-3}$ &  $4.6111\cdot 10^{-3}$\\
2         & $(693,543)$ & $3.2082 \cdot 10^{-3}$ &  $3.2136 \cdot 10^{-3}$\\
\end{tabular}
\caption{Example \ref{ex6}. Indices chosen by function eKG2 and global 
$K$-efficiency of both  ${\mathcal{G}}$ and ${\tilde{\mathcal{G}}}$, 
for $K=2,3,4,5$ and for $K=12:10:62$.} 
\label{T_2}
\end{table}
\end{example}

\section{Concluding remarks} \label{s8}
The adjacency matrix of a graph is a well-known tool for studying properties of a network 
defined by the graph. The path length matrix associated with a graph also sheds light on
properties of the network, but so far has not received much attention. A review of 
measures that can be defined in terms of the path length matrix is provided, and new such
measures are introduced. The sensitivity of the transmission of information to 
perturbations of the entries of the adjacency matrix is investigated. 

\section*{Data availability statement} 
Data sharing not is applicable to this article as no new datasets were generated during the 
current study.

\section*{Conflict of interest} 
The authors declare that they have no conflict of interest. 

\section*{Funding}
There is no funding available.

\section*{Authors' contributions}
All author contributed equally to the paper.

\bibliographystyle{plain}

\begin{thebibliography}{99}

\bibitem{air500}
{\em Dynamic Connectome Lab - Data Sets}. 
\url{https://sites.google.com/view/dynamicconnectomelab}

\bibitem{BBV}
A. Barrat, M. Barthelemy, and A. Vespignani, {\em Dynamical processes on complex networks},
Cambridge University Press, Oxford, 2008.

\bibitem{CH}
J. Clark and D. A. Holton, {\em A First Look at Graph Theory}, World Scientific, 
Singapore, 1991.


\bibitem{djnr_new}
O. De la Cruz Cabrera, J. Jin, S. Noschese, and L. Reichel, {\em Communication in complex 
networks}, Appl. Numer. Math., 172 (2022), pp. 186--205.

\bibitem{DMR}
O. De la Cruz Cabrera, M. Matar, and L. Reichel, {\em Analysis of directed networks via 
the matrix exponential}, J. Comput. Appl. Math., 355 (2019), pp. 182--192.

\bibitem{DMR2}
O. De la Cruz Cabrera, M. Matar, and L. Reichel, {\em Edge importance in a network via 
line graphs and the matrix exponential}, Numer. Algorithms, 83 (2020), pp. 807--832.

\bibitem{sm1} 
S. El-Halouy, S. Noschese, and L. Reichel, {\em Perron communicability and sensitivity of 
multilayer networks}, Numer. Algorithms, 92 (2023), pp. 597--617.

\bibitem{estrada2011structure}
E. Estrada, {\em The Structure of Complex Networks: Theory and Applications}, Oxford 
University Press, Oxford, 2011.

\bibitem{estrada2005subgraph}
E. Estrada and J. A. Rodriguez-Velazquez, {\em Subgraph centrality in complex networks}, 
Phys. Rev. E, 71 (2005), Art. 056103.

\bibitem{FRR}
C. Fenu, L. Reichel, and G. Rodriguez, {\em SoftNet: A package for the analysis of complex
networks}, Algorithms, 15 (2022), Art. 296.

\bibitem{GP}
R. L. Graham and H. O. Pollak, {\em On the addressing problem for loop switching}, 
The Bell System Technical Journal, 50 (1971), pp. 2495--2519.

\bibitem{L}
G. L. Litvinov, {\em Maslov dequantization, idempotent and tropical mathematics: A brief 
introduction}, J. Math. Sci., 140 (2007), pp. 426--444.

\bibitem{newman2004analysis}
M. E. J. Newman, {\em Analysis of weighted networks}, Phys. Rev E, 70 (2004), Art. 056131.

\bibitem{newman2010}
M. E. J. Newman, {\em Networks: An Introduction}, Oxford University Press, Oxford, 2010.

\bibitem{NP06}
S. Noschese and L. Pasquini, {\em Eigenvalue condition numbers: Zero-structured versus
traditional}, J. Comput. Appl. Math., 185 (2006), pp. 174--189.

\bibitem{NR}
S. Noschese and  L. Reichel, {\em Estimating and increasing the structural robustness of 
a network}, Numer. Linear Algebra Appl., 29 (2022), Art. e2418.

\bibitem{PNTM}
D. Plavsic, S. Nikolic, N. Trinajstic, and Z. Mihalic, {\em On the Harary index for the 
characterization of chemical graphs}, J. Math. Chem., 12 (1993), pp. 235--250.

\bibitem{R}
A. Ruhe, {\em The two-sided Arnoldi algorithm for nonsymmetric eigenvalue problems}, 
Matrix Pencils, eds. B. K\aa gstr\"om and A. Ruhe, Springer, Berlin, 1983, pp. 104--120.

\bibitem{taylor2009controllable}
A. Taylor and D.J. Higham, {\em A Controllable Test Matrix Toolbox for MATLAB}, ACM 
Trans. Math. Software, 35 (2009), pp. 26:1-26:17.

\bibitem{W}
J. H. Wilkinson, {\em Sensitivity of eigenvalues II}, Util. Math., 30 (1986), pp. 243--286.

\bibitem{ZH}
I. N. Zwaan and M. E. Hochstenbach, {\em Krylov-Schur-type restarts for the two-sided
Arnoldi method}, SIAM J. Matrix Anal. Appl., 38 (2017), pp. 297--321.

\end{thebibliography}

\end{document}